\RequirePackage{ifpdf}
\ifpdf 
\documentclass[pdftex]{sigma}
\else
\documentclass{sigma}
\fi

\renewcommand{\a}{\alpha}
\renewcommand{\b}{\beta}

\begin{document}

\allowdisplaybreaks

\renewcommand{\PaperNumber}{065}

\FirstPageHeading

\renewcommand{\thefootnote}{$\star$}

\ShortArticleName{The Rahman Polynomials Are Bispectral}

\ArticleName{The Rahman Polynomials Are Bispectral\footnote{This paper is a
contribution to the Vadim Kuznetsov Memorial Issue `Integrable
Systems and Related Topics'. The full collection is available at
\href{http://www.emis.de/journals/SIGMA/kuznetsov.html}{http://www.emis.de/journals/SIGMA/kuznetsov.html}}}

\Author{F. Alberto GR\"UNBAUM}

\AuthorNameForHeading{F.A. Gr\"unbaum}

\Address{Department of Mathematics, University of California, Berkeley,
CA 94720, USA}

\Email{\href{mailto:grunbaum@math.berkeley.edu}{grunbaum@math.berkeley.edu}}

\ArticleDates{Received February 01, 2007, in f\/inal form April
22, 2007; Published online May 03, 2007}

\Abstract{In a very  recent paper, M.~Rahman introduced a remarkable family of polynomials
in two variables as the eigenfunctions of the transition matrix for a nontrivial
Markov chain due to M.~Hoare and M.~Rahman.
I indicate here that these polynomials are bispectral. This should be just one of the many
remarkable properties enjoyed by these polynomials. For several challenges, including
f\/inding a general proof of some of the facts displayed here the reader should look
at the last section of this paper.}


\Classification{33C45; 22E45}

\begin{flushright}
\it Dedicated to the memory of Vadim Kutznetsov
\end{flushright}

\noindent
{\sl I first met Vadim at a meeting in Esterel, Canada back in 1994. I only saw him a few more times after that,
in Amsterdam, in Leeds, in Louvain-la-Neuve, and the last time in Leganes. Although we never wrote
any papers together we shared a number of common interests. I always found Vadim extremely helpful and
willing to share his ideas and his vast knowledge. He was also willing to put up with people whose work
was not in the mainstream and he usually would have some useful remark or a suggestion to make.
I like to think that he would have
enjoyed seeing this paper which builds on some remarkable work of another common friend I met for the first
time at the same Esterel meeting.}

\section{Preliminaries}

Starting about thirty years ago, M.~Hoare and M.~Rahman, published a number of 
papers studying a class of statistical models that were nicely associated
with ``classical'' orthogonal polynomials in one
variable. The reader can consult \cite{HR1,HR2,HR3,CHR}.

In a very recent paper \cite{HR4}, which I have seen before publication thanks to the 
kindness of M.~Rahman, they take a huge step and consider
the case of several variables. I propose to call these statistical models with the name Hoare--Rahman,
and to call the multivariable polynomials that arise with the name Rahman polynomials. I am thankful 
to Prof.~Rahman for several inspiring conversations on the contents of \cite{HR4}, 
in dif\/ferent cof\/fee houses and restaurants in Fremont, California on the last
days of 2006.

The model of interest here captures very well the idea of a game of chance where the player
takes a risk that can never improve his/her winnings but is then
given a second opportunity to improve his/her lot. I am convinced that models of this kind should be
of great interest in several areas of applied mathematics not only in physics (where they were conceived) but
also in biology, population dynamics, evolutionary models, etc. On top of their potential 
for applications the polynomials discovered by M.~Rahman are a
piece of beauty. I hope that this paper will incite other
people to look around for the gems contained in \cite{HR4}. I only recall the points from this paper that are needed
in the present context, urging the reader to go back to the source for inspiration.

For a systematic study of polynomials in several variables one should consult \cite{DX} as well as
the work on Macdonald polynomials of various kinds, see~\cite{M}. A look at the pioneering work of 
Tom Koornwinder, see for instance \cite{Ko}, is always
a very good idea.

\section{The Rahman polynomials}

Given generic values of the parameters $p_1$, $p_2$, $p_3$, $p_4$ (these could even be complex numbers) def\/ine
$t$, $u$, $v$, $w$ by means of
\begin{gather*}
t = \frac {(p_1+p_2)(p_1+p_3)}{p_1(p_1+p_2+p_3+p_4)}, \qquad u =
\frac {(p_1+p_3)(p_4+p_3)}{p_3(p_1+p_2+p_3+p_4)}, \\
v = \frac {(p_1+p_2)(p_2+p_4)}{p_2(p_1+p_2+p_3+p_4)}, \qquad w =
\frac {(p_4+p_2)(p_4+p_3)}{p_4(p_1+p_2+p_3+p_4)}.
\end{gather*}

Now we can def\/ine with M.~Rahman the polynomials in $x$, $y$ given by
\[
P_{m,n}(x,y) = \sum_i \sum_j \sum_k \sum_{\ell} \frac
{(-m)_{i+j}(-n)_{k+\ell}(-x)_{i+k}(-y)_{j+\ell}}{i!j!k!\ell!(-N)_{i+j+k+\ell}}
\, t^iu^jv^kw^{\ell}.
\]
I will not try to explain M.~Rahman's path to these polynomials, but the reader should consult~\cite{HR4}. 
In particular one f\/inds in~\cite{HR4} a discussion of the connection with the $9-j$ symbols used in quantum mechanics.

\section[The Hoare-Rahman model and the kernel ${\mathcal K}$]{The Hoare--Rahman model and the kernel 
$\boldsymbol{\mathcal K}$}

    Denote by 
\[
    b(k,N,\alpha)={N \choose k} \alpha^k (1-\alpha)^{N-k}
\]
the usual binomial distribution giving the
    probability of $k$ successes in $N$ repetitions of an experiment with two possible outcomes: success
    (with probability $\alpha$) and failure (with probability $1-\alpha$). Likewise denote by
$b_2(i_1,i_2,N;\beta_1,\beta_2)$ the trinomial distribution with parameters $\beta_1$, $\beta_2$.

This means that that for $i_1$, $i_2$ nonnegative and such that $i_1+i_2 \leq N$ we have
\[
b_2(i_1,i_2,N;\beta_1,\beta_2) = \beta_1^{i_1}\beta_2^{i_2}(1-\beta_1-\beta_2)^{N-i_1-i_2} \left( \frac
{N!}{i_1!i_2!(N-i_1-i_2)!} \right).
\]

    Def\/ine now a transition probability kernel 
    \[
    {\mathcal K}(j_1,j_2;i_1,i_2)
    \] in the space of
    non-negative integer physical
    variables $i_1$, $i_2$ and $j_1$, $j_2$ such that
\[
i_1+i_2 \leq N, \qquad j_1+j_2 \leq N
\]
by means of the expression
\begin{gather*}
{\mathcal K}(j_1,j_2;i_1,i_2) = \sum_{k_1=0}^{\min(i_1,j_1)}
\sum_{k_2=0}^{\min(i_2,j_2)}
b(k_1,i_1;\alpha_1)b(k_2,i_2;\alpha_2)\\
\phantom{{\mathcal K}(j_1,j_2;i_1,i_2) =}{}\times b_2(j_1-k_1,j_2-k_2,N-k_1-k_2;\beta_1,\beta_2).
\end{gather*}

The expression above is written so as to agree with the one in~\cite{HR4}. Notice that its interpretation is that of the probability of a transition from
the state $(i_1,i_2)$ to the new
sate $(j_1,j_2)$ in one step (consisting of three tosses) of our process, as will be explained in the next section.

    At this point the parameters $\alpha_1$, $\alpha_2$, $\beta_1$, $\beta_2$ are free to lie anywhere in the interval
    $[0,1]$, with the sum of $\beta_1$ and $\beta_2$ not exceeding the value $1$.

    The corresponding Markov chain is not a ``nearest neighbour one'' in the sense that the classical
    Ehrenfest urn model, see \cite{E,F,K,SchK}, or all other classical models are. Here the transition
    probability matrix has {\em all} its entries positive and therefore there is no sense in which a
    one-step transition
    moves you only to your immediate neighbours.  I f\/ind it remarkable that the complete spectral problem
    can be solved for such a kernel, but this is exactly what the Rahman polynomials do, at least under certain
    extra conditions.

\section{Poker dice}

A concrete example of the model is provided by a gambling game known as Poker dice.
The version considered by M.~Hoare and M.~Rahman goes as follows.
You have $N$ identical dice with $n$ faces each. There are three kinds of faces:
red, black and ``blanks''. The probabilities of getting a red or a black face when
rolling one dice are $\alpha_1$ and $\alpha_2$, respectively. The state of the system is given
by $(i_1,i_2)$, the number of red and black faces among the $N$ dice. Think of the remaining $N-i_1-i_2$
dice as failures.

Start form the state $(i_1,i_2)$ and get ready to roll three times. The resulting
number of reds and blacks after these three rools will be denoted by $(j_1,j_2)$. 
This will constitute one time step in the Markov chain to be described
now.

The $i_1$ dice that
show red are tossed and yield $k_1$ red faces. The $i_2$ dice that showed a
black face are tossed and yield $k_2$ blacks. Now you take the initial
$N-i_1-i_2$ blanks and combine them with the $i_1-k_1$ dice that failed to reproduce
their original red color and
the $i_2-k_2$ dice that failed to reproduce their black color.
The total number of ``failures'' is now $N-k_1-k_2$.
These dice are thrown by themselves: they yield a certain number of red  and black faces which are
added to the $k_1$ red faced dice and $k_2$ black faced dice that were saved earlier. This total number of
reds is now $j_1$ and the total number of blacks is $j_2$. To obtain full generality
the probabilities of getting a red or a black for each individual dice on this last
toss can be taken to be $\beta_1$ and $\beta_2$, not necessarily the old values
$\alpha_1$ and $\alpha_2$.

The case of $N=5$ is referred to as standard ``poker dice'' in \cite{HR4} and will
be used to illustrate later on the main point in this paper, namely the
bispectral property of the Rahman polynomials. We have verif\/ied this property for many other values of $N$, 
but we do not have a nice general expression
for the operator in question, or a nice proof of its existence for general $N$.

\section{The orthogonality weight}

In \cite{HR4} one sees that the Rahman polynomials are orthogonal with respect to the trinomial
distribution 
\[
b_2(x,y;N,\eta_1,\eta_2),
\] where the parameters $\eta_1$, $\eta_2$ are given in terms of
$p_1$, $p_2$, $p_3$, $p_4$ by means of
\begin{gather*}
\eta_1 = \frac {p_1p_2(p_1+p_2+p_3+p_4)}{(p_1+p_2)(p_1+p_3)(p_2+p_4)}, \qquad
\eta_2 = \frac {p_3p_4(p_1+p_2+p_3+p_4)}{(p_1+p_3)(p_4+p_2)(p_4+p_3)}.
\end{gather*}

Given this weight, for an arbitrary choice of the parameters  $\eta_1$, $\eta_2$, 
one could blindly decide to ``order'' the set of monomials 
\[
x^i y^j
\]
 and
then apply
the Gram--Schmidt  process to obtain orthogonal polynomials.
Of course there is no unique way to order these monomials and each ordering will produce a dif\/ferent
family of polynomials.

One of the features of the work
of M.~Rahman is that it automatically gives one set of polynomials that are already orthogonal and,
as we see in \cite{HR4} and here, enjoy remarkable properties. Moreover he produces a four parameter family
of these polynomials, having started from a the two parameters $\eta_1$, $\eta_2$. There is clearly a lot of
inspiration behind this construction, and the reader should consult~\cite{HR4}.

\section{The bispectral property}

Recall that ${\mathcal K}$ depends on one set of four parameters and the Rahman
polynomials on a dif\/ferent set of four parameters.

  As shown in \cite{HR4}, the Rahman polynomials satisfy the eigenvalue relation
\[
{\mathcal K} P_{m,n} (x,y) = \lambda_{m,n} P_{m,n}(x,y)
\]
for $(x,y)$ in ``physical space'', i.e.\ $x$, $y$ are integers satisfying  $x \geq 0$, $y \geq 0$, $x+y \leq N$,
with
\begin{gather*}
\lambda_{m,n} = \left\{ \frac {1-\alpha_1\eta_1 - \alpha_2\eta_2 -
\eta_1t(1-\alpha_1) - \eta_2u(1-\alpha_2)}{1-\alpha_1\eta_1 -
\alpha_2\eta_2} \right\}^m \\
\phantom{\lambda_{m,n} =}{}\times \left\{ \frac
{1-\alpha_1\eta_1-\alpha_2\eta_2-\eta_1v(1-\alpha_1)-\eta_2w(1-\alpha_2)}{1-\alpha_1\eta_1-\alpha_2\eta_2}
\right\}^n \\
\phantom{\lambda_{m,n}}{}= (1-\beta_1t-\beta_2u)^m(1-\beta_1v-\beta_2w)^n.
\end{gather*}

More explicitly, with the def\/initions made above this reads,
\[
\sum_{j_1,j_2}  K(j_1,j_2;i_1,i_2) P_{m,n}(j_1,j_2) = \lambda_{m,n} P_{m,n}(i_1,i_2).
\]
  For this to be true we need to impose some relations connecting the parameters $p_1$,
  $p_2$, $p_3$, $p_4$ and
  $\alpha_1$, $\alpha_2$, $\beta_1$, $\beta_2$. This issue has been addressed in~\cite{HR4}.

  In the considerations above
  we are thinking of these polynomials as functions of the physical
  space variables $x$, $y$ for f\/ixed values
  of the frequency parameters $m$, $n$.

  In the spirit of \cite{DG} one could ask if they are the eigenvectors of another spectral problem
  when considered, for f\/ixed $x$, $y$, as functions of $m$, $n$. If this were true we would have an
  analog of the three-term recursion relation satisf\/ied by orthogonal polynomials of one variable.
  In the continuous-continuous case this question was raised in~\cite{DG} for a second order Schr\"odinger
  type dif\/ferential operator. For many ramif\/ications of this problem the reader can see~\cite{HK}.

  We give below the simplest recursion that we have been able to f\/ind. To be very concrete we 
  dis\-play the result in the case of $N=5$ which corresponds to the Poker dice example described above.
  The operator in question consists of multiplication by the function $(p_1{+}p_2)x{-}(p_3{+}p_4)y$
  and it gives a ``f\/ive-point'' formula, where most points in the ``spectral space'' have {\em four}
  neighbours. There are, as indicated above, similar recursions for other values of $N$ but we do not have a general expression that would cover all cases.

  In the case of $N=5$ the frequency space (and also the physical space) looks as follows

\centerline{\includegraphics[width=6cm]{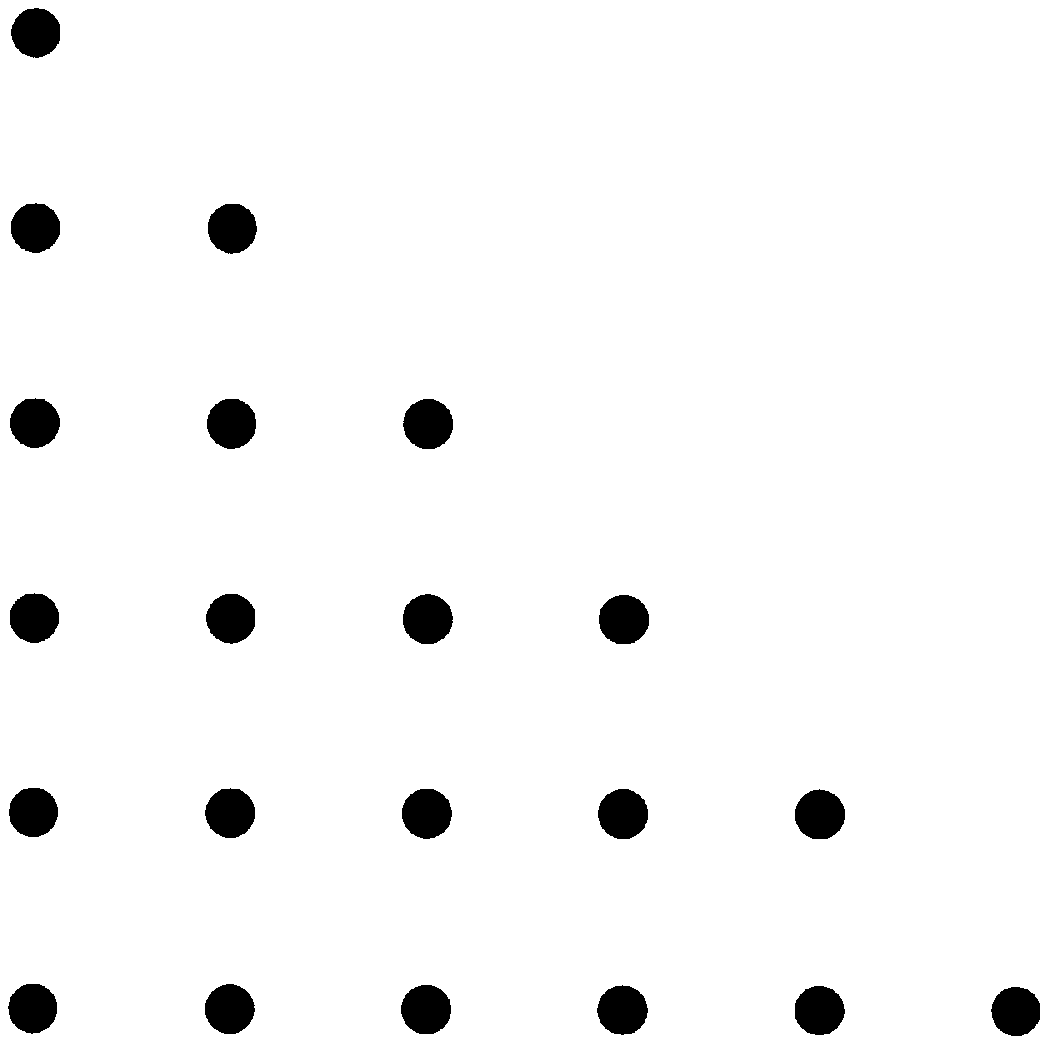}}

\medskip

  If we declare two points to be adjacent when they can be obtained by going one step either in the
  vertical or horizontal directions, and we order the $21$ states as follows
\begin{gather*}
(0,0),(1,0),(2,0),(3,0),(4,0),(5,0),(0,1),(1,1),(2,1),(3,1),(4,1),\\
(0,2),(1,2),(2,2),(3,2),(0,3),(1,3),(2,3),(0,4),(1,4),(0,5)
\end{gather*}
then one can work out easily the adjacency matrix.

The matrix ${\mathcal B}$ to be detailed below will have nonzero entries exactly where the adjacency matrix has a $1$. It is therefore as much of a {\em
local} matrix as one can get in
this setup. It will correspond to the analog of the {\em three term recursion
relation} of the usual situation when one deals with orthogonal polynomials in
one variable.

We claim that the Rahman polynomials, as functions of $m$, $n$ satisfy the dif\/ference equation
\[
{\mathcal B} P_{m,n} (x,y) = ((p_1+p_2) x - (p_3+p_4) y) P_{m,n}(x,y)
\]
as long as $(x,y)$ lie in physical space.

Here the matrix ${\mathcal B}$ is made up of rows as specif\/ied one at a time below.
 Notice that in each row there is an element denoted by $s_i$  for which we do not
 give an expression. This is not needed since the sum of all the entries in each row
 equals the value zero.
 \[
[s_1,a_1,0,0,0,0,b_1,0,0,0,0,0,0,0,0,0,0,0,0,0,0]
\]
with
\begin{gather*}
a_1 = \frac
{5p_1p_3(p_4+p_2)(p_4+p_3+p_2+p_1)}{(p_3+p_1)(p_1p_4-p_2p_3)}, \qquad
b_1 = -\frac {5p_2(p_3+p_1)p_4(p_4+p_3+p_2+p_1)}{(p_4+p_2)(p_1p_4-p_2p_3)};
\\
[a_2,s_2,c_2,0,0,0,0,d_2,0,0,0,0,0,0,0,0,0,0,0,0,0]
\end{gather*}
with
\begin{gather*}
a_2 = \frac {p_1p_4-p_2p_3}{p_3+p_1}, \qquad
c_2 = \frac
{4p_1p_3(p_4+p_2)(p_4+p_3+p_2+p_1)}{(p_3+p_1)(p_1p_4-p_2p_3)}, \\
d_2 = -\frac {4p_2(p_3+p_1)p_4(p_4+p_3+p_2+p_1)}{(p_4+p_2)(p_1p_4-p_2p_3)};
\\
[0,a_3,s_3,c_3,0,0,0,0,d_3,0,0,0,0,0,0,0,0,0,0,0,0]
\end{gather*}
with
\begin{gather*}
a_3 = \frac {2(p_1p_4-p_2p_3)}{p_3+p_1}, \qquad
c_3 = \frac
{3p_1p_3(p_4+p_2)(p_4+p_3+p_2+p_1)}{(p_3+p_1)(p_1p_4-p_2p_3)}, \\
d_3 = -\frac {3p_2(p_3+p_1)p_4(p_4+p_3+p_2+p_1)}{(p_4+p_2)(p_1p_4-p_2p_3)};
\\
[0,0,a_4,s_4,c_4,0,0,0,0,d_4,0,0,0,0,0,0,0,0,0,0,0]
\end{gather*}
with
\begin{gather*}
a_4 = \frac {3(p_1p_4-p_2p_3)}{p_3+p_1}, \qquad
c_4 = \frac
{2p_1p_3(p_4+p_2)(p_4+p_3+p_2+p_1)}{(p_3+p_1)(p_1p_4-p_2p_3)}, \\
d_4 = -\frac {2p_2(p_3+p_1)p_4(p_4+p_3+p_2+p_1)}{(p_4+p_2)(p_1p_4-p_2p_3)};
\\
[0,0,0,a_5,s_5,c_5,0,0,0,0,e_5,0,0,0,0,0,0,0,0,0,0]
\end{gather*}
with
\begin{gather*}
a_5 = \frac {4(p_1p_4-p_2p_3)}{p_3+p_1}, \qquad
c_5 = \frac
{p_1p_3(p_4+p_2)(p_4+p_3+p_2+p_1)}{(p_3+p_1)(p_1p_4-p_2p_3)}, \\
e_5 = -\frac {p_2(p_3+p_1)p_4(p_4+p_3+p_2+p_1)}{(p_4+p_2)(p_1p_4-p_2p_3)};
\\
[0,0,0,0,a_6,s_6,0,0,0,0,0,0,0,0,0,0,0,0,0,0,0]
\end{gather*}
with
\begin{gather*}
a_6 = \frac {5(p_1p_4-p_2p_3)}{p_3+p_1};
\\
[d_7,0,0,0,0,0,s_7,c_7,0,0,0,a_7,0,0,0,0,0,0,0,0,0]
\end{gather*}
with
\begin{gather*}
a_7 = -\frac
{4p_2(p_3+p_1)p_4(p_4+p_3+p_2+p_1)}{(p_4+p_2)(p_1p_4-p_2p_3)}, \\
c_7 = \frac
{4p_1p_3(p_4+p_2)(p_4+p_3+p_2+p_1)}{(p_3+p_1)(p_1p_4-p_2p_3)}, \qquad
d_7 = -\frac {p_1p_4-p_2p_3}{p_4+p_2};
\\
[0,a_8,0,0,0,0,b_8,s_8,d_8,0,0,0,e_8,0,0,0,0,0,0,0,0]
\end{gather*}
with
\begin{gather*}
a_8 = -\frac {p_1p_4-p_2p_3}{p_4+p_2}, \qquad
b_8 = \frac {p_1p_4-p_2p_3}{p_3+p_1}, \\
d_8 = \frac
{3p_1p_3(p_4+p_2)(p_4+p_3+p_2+p_1)}{(p_3+p_1)(p_1p_4-p_2p_3)}, \qquad
e_8 = -\frac {3p_2(p_3+p_1)p_4(p_4+p_3+p_2+p_1)}{(p_4+p_2)(p_1p_4-p_2p_3)};
\\
[0,0,e_9,0,0,0,0,c_9,s_9,a_9,0,0,0,d_9,0,0,0,0,0,0,0]
\end{gather*}
with
\begin{gather*}
a_9 = \frac
{2p_1p_3(p_4+p_2)(p_4+p_3+p_2+p_1)}{(p_3+p_1)(p_1p_4-p_2p_3)}, \qquad
c_9 = \frac {2(p_1p_4-p_2p_3)}{p_3+p_1}, \\
d_9 = -\frac
{2p_2(p_3+p_1)p_4(p_4+p_3+p_2+p_1)}{(p_4+p_2)(p_1p_4-p_2p_3)}, \qquad
e_9 = -\frac {p_1p_4-p_2p_3}{p_4+p_2};
\\
[0,0,0,d_{10},0,0,0,0,a_{10},s_{10},c_{10},0,0,0,e_{10},0,0,0,0,0,0]
\end{gather*}
with
\begin{gather*}
a_{10} = \frac {3(p_1p_4-p_2p_3)}{p_3+p_1}, \qquad
c_{10} = \frac
{p_1p_3(p_4+p_2)(p_4+p_3+p_2+p_1)}{(p_3+p_1)(p_1p_4-p_2p_3)}, \\
d_{10} = -\frac {p_1p_4 - p_2p_3}{p_4+p_2}, \qquad
e_{10} = -\frac
{p_2(p_3+p_1)p_4(p_4+p_3+p_2+p_1)}{(p_4+p_2)(p_1p_4-p_2p_3)};
\\
[0,0,0,0,b_{11},0,0,0,0,a_{11},s_{11},0,0,0,0,0,0,0,0,0,0]
\end{gather*}
with
\begin{gather*}
a_{11} = \frac {4(p_1p_4-p_2p_3)}{p_3+p_1}, \qquad
b_{11} = -\frac {p_1p_4 - p_2p_3}{p_4+p_2};
\\
[0,0,0,0,0,0,a_{12},0,0,0,0,s_{12},d_{12},0,0,c_{12},0,0,0,0,0]
\end{gather*}
with
\begin{gather*}
a_{12} = -\frac {2(p_1p_4-p_2p_3)}{p_4+p_2}, \qquad
c_{12} = -\frac
{3p_2(p_3+p_1)p_4(p_4+p_3+p_2+p_1)}{(p_4+p_2)(p_1p_4-p_2p_3)}, \\
d_{12} = \frac
{3p_1p_3(p_4+p_2)(p_4+p_3+p_2+p_1)}{(p_3+p_1)(p_1p_4-p_2p_3)};
\\
[0,0,0,0,0,0,0,c_{13},0,0,0,a_{13},s_{13},d_{13},0,0,e_{13},0,0,0,0]
\end{gather*}
with
\begin{gather*}
a_{13} = \frac {p_1p_4-p_2p_3}{p_3+p_1}, \qquad
c_{13} = -\frac {2(p_1p_4-p_2p_3)}{p_4+p_2}, \\
d_{13} = \frac
{2p_1p_3(p_4+p_2)(p_4+p_3+p_2+p_1)}{(p_3+p_1)(p_1p_4-p_2p_3)}, \qquad
e_{13} = -\frac
{2p_2(p_3+p_1)p_4(p_4+p_3+p_2+p_1)}{(p_4+p_2)(p_1p_4-p_2p_3)};
\\
[0,0,0,0,0,0,0,0,d_{14},0,0,0,a_{14},s_{14},c_{14},0,0,e_{14},0,0,0]
\end{gather*}
with
\begin{gather*}
a_{14} = \frac {2(p_1p_4-p_2p_3)}{p_3+p_1}, \qquad
c_{14} = \frac
{p_1p_3(p_4+p_2)(p_4+p_3+p_2+p_1)}{(p_3+p_1)(p_1p_4-p_2p_3)}, \\
d_{14} = -\frac {2(p_1p_4-p_2p_3)}{p_4+p_2}, \qquad
e_{14} = -\frac
{p_2(p_3+p_1)p_4(p_4+p_3+p_2+p_1)}{(p_4+p_2)(p_1p_4-p_2p_3)};
\\
[0,0,0,0,0,0,0,0,0,a_{15},0,0,0,c_{15},s_{15},0,0,0,0,0,0]
\end{gather*}
with
\begin{gather*}
a_{15} = -\frac {2(p_1p_4-p_2p_3)}{p_4+p_2}, \qquad
c_{15} = \frac {3(p_1p_4-p_2p_3)}{p_3+p_1};
\\
[0,0,0,0,0,0,0,0,0,0,0,a_{16},0,0,0,s_{16},c_{16},0,d_{16},0,0]
\end{gather*}
with
\begin{gather*}
a_{16} = -\frac {3(p_1p_4-p_2p_3)}{p_4+p_2}, \qquad
c_{16} = \frac
{2p_1p_3(p_4+p_2)(p_4+p_3+p_2+p_1)}{(p_3+p_1)(p_1p_4-p_2p_3)}, \\
d_{16} = -\frac
{2p_2(p_3+p_1)p_4(p_4+p_3+p_2+p_1)}{(p_4+p_2)(p_1p_4-p_2p_3)};
\\
[0,0,0,0,0,0,0,0,0,0,0,0,e_{17},0,0,a_{17},s_{17},c_{17},0,d_{17},0]
\end{gather*}
with
\begin{gather*}
a_{17} = \frac {p_1p_4-p_2p_3}{p_3+p_1}, \qquad
c_{17} = \frac
{p_1p_3(p_4+p_2)(p_4+p_3+p_2+p_1)}{(p_3+p_1)(p_1p_4-p_2p_3)}, \\
d_{17} = -\frac
{p_2(p_3+p_1)p_4(p_4+p_3+p_2+p_1)}{(p_4+p_2)(p_1p_4-p_2p_3)}, \qquad
e_{17} = -\frac {3(p_1p_4-p_2p_3)}{p_4+p_2};
\\
[0,0,0,0,0,0,0,0,0,0,0,0,0,a_{18},0,0,c_{18},s_{18},0,0,0]
\end{gather*}
with
\begin{gather*}
a_{18} = -\frac {3(p_1p_4-p_2p_3)}{p_4+p_2}, \qquad
c_{18} = \frac {2(p_1p_4-p_2p_3)}{p_3+p_1};
\\
[0,0,0,0,0,0,0,0,0,0,0,0,0,0,0,c_{19},0,0,s_{19},a_{19},d_{19}]
\end{gather*}
with
\begin{gather*}
a_{19} = \frac
{p_1p_3(p_4+p_2)(p_4+p_3+p_2+p_1)}{(p_3+p_1)(p_1p_4-p_2p_3)}, \\
c_{19} = -\frac {4(p_1p_4-p_2p_3)}{p_4+p_2}, \qquad
d_{19} = -\frac
{p_2(p_3+p_1)p_4(p_4+p_3+p_2+p_1)}{(p_4+p_2)(p_1p_4-p_2p_3)};
\\
[0,0,0,0,0,0,0,0,0,0,0,0,0,0,0,0,a_{20},0,c_{20},s_{20},0]
\end{gather*}
with
\begin{gather*}
a_{20} = -\frac {4(p_1p_4-p_2p_3)}{p_4+p_2}, \qquad
c_{20} = \frac {p_1p_4-p_2p_3}{p_3+p_1};
\\
[0,0,0,0,0,0,0,0,0,0,0,0,0,0,0,0,0,0,a_{21},0,s_{21}]
\end{gather*}
with
\[
a_{21} = -\frac {5(p_1p_4-p_2p_3)}{p_4+p_2}.
\]

  We note that these relations hold for arbitrary (generic) values of $p_1$, $p_2$, $p_3$, $p_4$.

  If one considers multiplication by $x$ or by $y$ we get a ``seven-point'' formula. Of course it would be
  interesting, in the spirit of \cite{CG} to consider the full algebra of dif\/ference operators going
  along with the Rahman polynomials.

  Thinking along these lines one could also consider the full algebra of dif\/ference operators that go along
  with the Rahman polynomials when thought of as functions of $x$, $y$. This simply means the algebra of
  all matrices that commute with the kernel ${\mathcal K}$. While the original kernel ${\mathcal K}$ has
  a very nice probabilistic interpretation there is no reason why other operators in this algebra should
  allow such an interpretation. More on this in the next section.

\section{Global vs. local}

   One should naturally think of the kernel ${\mathcal K}$ given earlier as a global
   operator, one that in this discrete case is represented by a {\em
full matrix} and that in the continuous case would be represented by an integral operator.
At the core of many considerations leading up to the bispectral problem is the
issue of the existence of a local operator that would commute with the global
one. This remarkable piece of work is forever associated with the names of D.~Slepian, 
H.~Landau and H.~Pollack. This is discussed for instance in \cite{G2,G3}.

Here, as an illustration we show that this question has a positive answer in the
setup of Hoare and Rahman. Take $N=3$ in which case we can see, for arbitrary
choice of the free parameters $\alpha_1$, $\alpha_2$, $\beta_1$, $\beta_2$ that the
matrix given below is essentially the only one that commutes with the kernel
${\mathcal K}$. Once again, we have checked that this holds for other values of $N$, and
this case is used here only for illustration purposes.

In the spirit of the work just alluded to one should consider ${\mathcal K}$ as the analog 
of the operator of time-band-time
limiting and the operator below as the analog of the second order dif\/ferential operator that comes from separation of
variables for the Laplacian in ${\mathbb R}^3$ in prolate spheroidal coordinates.

In agreement with the case $N=5$ considered as an illustration earlier, we order the 
states in physical space as follows
\[
(0,0),(1,0),(2,0),(3,0),(0,1),(1,1),(2,1),(0,2),(1,2),(0,3).
\]

Then, up to trivial scaling and shifts the only {\em local} matrix that
commutes with ${\mathcal K}$ is given by
\[
\begin{bmatrix}
x_1 & x_2 & 0 & 0 & x_3 & 0 & 0 & 0 & 0 & 0 \\
x_4 & x_5 & x_6 & 0 & 0 & x_7 & 0 & 0 & 0 & 0 \\
0 & x_8 & x_9 & x_{10} & 0 & 0 & x_{11} & 0 & 0 & 0 \\
0 & 0 & x_{12} & x_{13} & 0 & 0 & 0 & 0 & 0 & 0 \\
x_{14} & 0 & 0 & 0 & x_{15} & x_{16} & 0 & x_{17} & 0 & 0 \\
0 & x_{18} & 0 & 0 & x_{19} & x_{20} & x_{21} & 0 & x_{22} & 0 \\
0 & 0 & x_{23} & 0 & 0 & x_{24} & x_{25} & 0 & 0 & 0 \\
0 & 0 & 0 & 0 & x_{26} & 0 & 0 & x_{27} & x_{28} & x_{29} \\
0 & 0 & 0 & 0 & 0 & x_{30} & 0 & x_{31} & x_{32} & 0 \\
0 & 0 & 0 & 0 & 0 & 0 & 0 & x_{33} & 0 & x_{34}
\end{bmatrix}
\]
as long as we choose its entries as follows
\begin{gather*}
x_1 = 3(\a_1\a_2\b_2 - 2\a_1\b_2 + \a_1\a_2\b_1 - \a_2\b_1 - \a_1\b_1
- \a_1\a_2 + \a_1)/(\a_1\b_2),
\\
x_2 = \frac {3\a_2\b_1}{\a_1\b_2}
\qquad 
x_3 = 3,
\qquad
x_4 = \frac {(\a_1 - 1)\a_2(\b_2 + \b_1 - 1)}{\a_1\b_2},
\\
x_5 = (2\a_1\a_2\b_2 + \a_2\b_2 - 5\a_1\b_2 + 2\a_1\a_2\b_1 -
\a_2\b_1 - 3\a_1\b_1 - 2\a_1\a_2 - \a_2 + 3\a_1)/(\a_1\b_2),
\\
x_6 = \frac {2\a_2\b_1}{\a_1\b_2},
\qquad
x_7 = 2,
\qquad
x_8 = \frac {2(\a_1 - 1)\a_2(\b_2 + \b_1 - 1)}{\a_1\b_2},
\\
x_9 = (\a_1\a_2\b_2 + 2\a_2\b_2 - 4\a_1\b_2 + \a_1\a_2\b_1 + \a_2\b_1
- 3\a_1\b_1 - \a_1\a_2 - 2\a_2 + 3\a_1)/(\a_1\b_2),
\\
x_{10} = \frac {\a_2\b_1}{\a_1\b_2},
\qquad
x_{11} = 1,
\qquad
x_{12} = \frac {3(\a_1 - 1)\a_2(\b_2 + \b_1 - 1)}{\a_1\b_2},
\\
x_{13} = \frac {3(\a_2 - \a_1)(\b_2 + \b_1 - 1)}{\a_1\b_2},
\qquad
x_{14} = \frac {(\a_2 - 1)(\b_2 + \b_1 - 1)}{\b_2},
\\
x_{15} = 2(\a_1\a_2\b_2 - 2\a_1\b_2 + \a_1\a_2\b_1 - \a_2\b_1 -
\a_1\b_1 - \a_1\a_2 + \a_1)/(\a_1\b_2),
\\
x_{16} = \frac {2\a_2\b_1}{\a_1\b_2},
\qquad
x_{17} = 2,
\qquad
x_{18} = \frac {(a_2 - 1)(\b_2 + b_1 - 1)}{\b_2},
\\
x_{19} = \frac {(\a_1 - 1)\a_2(\b_2 + \b_1 - 1)}{\a_1\b_2},
\\
x_{20} = (\a_1\a_2\b_2 + \a_2\b_2 - 3\a_1\b_2 + \a_1\a_2\b_1 -
2\a_1\b_1 - \a_1\a_2 - \a_2 + 2\a_1)/(\a_1\b_2),
\\
x_{21} = \frac {\a_2\b_1}{a_1\b_2},\qquad
x_{22} = 1,
\qquad
x_{23} = \frac {(\a_2 - 1)(\b_2 + \b_1 - 1)}{\b_2},
\\
x_{24} = \frac {2(\a_1 - 1)\a_2(\b_2 + \b_1 - 1)}{\a_1\b_2},
\qquad
x_{25} = \frac {2(\a_2 - \a_1)(\b_2 + \b_1 - 1)}{\a_1\b_2},
\\
x_{26} = \frac {2(\a_2 - 1)(\b_2 + \b_1 - 1)}{\b_2},
\\
x_{27} = (\a_1\a_2\b_2 - 2\a_1\b_2 + \a_1\a_2\b_1 - \a_2\b_1 -
\a_1\b_1 - \a_1\a_2 + \a_1)/(\a_1\b_2),
\\
x_{28} = \frac {\a_2\b_1}{\a_1\b_2},
\qquad
x_{29} = 1,
\qquad
x_{30} = \frac {2(\a_2 - 1)(\b_2 + \b_1 - 1)}{\b_2},
\\
x_{31} = \frac {(\a_1 - 1)\a_2(\b_2 + \b_1 - 1)}{\a_1\b_2},
\qquad
x_{32} = \frac {(\a_2 - \a_1)(\b_2 + \b_1 - 1)}{\a_1\b_2},
\\
x_{33} = \frac {3(\a_2 - 1)(\b_2 + \b_1 - 1)}{\b_2},
\qquad
x_{34} = 0.
\end{gather*}

\section{Some closing challenges}

  Here is a challenge that is motivated by recent work on matrix valued orthogonal polynomials.
  This subject was initiated by M.G.~Krein a long time ago, see \cite{K1,K2}. There has been quite 
  a~bit of work in terms of producing explicit examples of these polynomials.
  For a representative sample see \cite{CG,DG1,DvA,G,GPT1,GPT2,GPT3,SvA,T2}. Most of this work deals with
  looking for examples that would enjoy the appropriate bispectral property.
  On the surface of it, there is little relation between these polynomials and the topics discussed here.
  But appearances can be deceiving.
  In a recent paper~\cite{G1} I noticed how the example in~\cite{G} has a probabilistic interpretation
  (still waiting to be unravelled) and leads to nonnearest neighbours transitions on the nonnegative
  integers. This goes beyond the usual birth-and-death processes for which there is a well known relation
  to orthogonal polynomials, see~\cite{AAR,AW,I,ILMV,KMcG,KMcG2}. In the same paper \cite{G1} one f\/inds examples
  arising in the theory of matrix valued orthogonal polynomials (in one variable) that give rise to
  Markov chains in two dimensional discrete space.

  Another challenge would be to f\/ind the appropriate concrete algebraic geometry that lies behind the relation among
  the variables $\alpha_1$, $\alpha_2$, $\beta_1$, $\beta_2$ and $p_1$, $p_2$, $p_3$, $p_4$. In particular
  this could have some bearing on the situation discussed in \cite{HR4} where the spectrum  of
  ${\mathcal K}$ is degenerate, see~(5.15).

  One challenge mentioned already a couple of times consists in writing down the general
  formula for a f\/ive-term recursion formula that is displayed here only for $N=5$.

\subsection*{Acknowledgments}

I am very thankful to a couple of referees who read the paper  with great care and pointed out typos as well as ways
  to improve the presentation.
  The author was supported in part by NSF Grant \# 0603901.

\pdfbookmark[1]{References}{ref}
\LastPageEnding

\end{document}